\documentclass[12pt]{elsart}
\usepackage{amssymb,amsmath,latexsym}

\newcommand{\proof}{\noindent\underline{Proof.}\ }
\newcommand{\zentr}{\mbox{{\sf Y}}}

\newcommand{\fa}{\mbox{ for}\mbox{ all }}
\newcommand{\ba}{\begin{array}}
\newcommand{\ea}{\end{array}}

\newcommand{\Z}{{\mathbb{Z}}}
\newcommand{\Q}{{\mathbb{Q}}}
\newcommand{\F}{{\mathbb{F}}}

\newcommand{\C}{{\mathbb{C}}}

\newcommand{\trace}{\tau }
\DeclareMathOperator{\wt}{wt}

\DeclareMathOperator{\Inv}{Inv}
\DeclareMathOperator{\cwe}{cwe}

\DeclareMathOperator{\GL}{GL}

\DeclareMathOperator{\SL}{SL}
\DeclareMathOperator{\Aut}{Aut}
\DeclareMathOperator{\QR}{QR}
\DeclareMathOperator{\id}{id}
\DeclareMathOperator{\Out}{Out}
\DeclareMathOperator{\Gal}{Gal}

\newcommand{\1}{{\bf 1}}
\newcommand{\eb}{\phantom{zzz}\hfill{$\square $}\smallskip}

\renewcommand{\em}{\sf}

\begin{document}
\begin{center}
{\Large {\bf  Complete Weight Enumerators of Generalized Doubly-Even Self-Dual Codes}} \\
\vspace{1.5\baselineskip}
{\em Gabriele Nebe}\footnote{
Abteilung Reine Mathematik, Universit\"at Ulm,
89069 Ulm, Germany,
nebe@mathematik.uni-ulm.de },
{\em H.-G. Quebbemann}\footnote{
Fachbereich Mathematik,
Universit\"at Oldenburg,
26111 Oldenburg, Germany,
quebbemann@mathematik.uni-oldenburg.de},
{\em E. M. Rains}\footnote{
Mathematics Department,
University of California Davis, 
Davis, CA 95616, USA,
rains@math.ucdavis.edu}
and {\em N. J. A. Sloane}\footnote{
Information Sciences Research, AT\&T Shannon Labs,
Florham Park, NJ 07932-0971, USA,
njas@research.att.com }
\end{center}

\vspace{1.5\baselineskip}

\noindent
{\bf ABSTRACT}
For any $q$ which is a power of 2
we describe a finite subgroup of $\GL_q(\C)$ under which the
complete weight enumerators of generalized doubly-even self-dual codes
over $\F_q$ are invariant. An explicit description of the invariant ring
and some applications to extremality of such codes are obtained in the
case $q=4$.
\vspace{.5\baselineskip}
\vspace{1\baselineskip}

\section{Introduction}
In 1970 Gleason \cite{Gle70} described a finite complex linear group of
degree $q$ under which the complete weight enumerators of self-dual codes
over $\F_q$ are invariant. While for odd $q$ this group is a
double or quadruple cover of $\SL_2(\F_q)$, for even $q \geq 4$ it is solvable of order
$4q^2(q-1)$ (compare \cite{klemm}). 
For even $q$ it is only when $q=2$
that the seemingly exceptional type of doubly-even self-dual binary codes
leads to a larger group.

In this paper we study a generalization of doubly-even codes to
the non-binary case which was introduced
in \cite{quebbemann}. A linear code of length $n$ over $\F_q$ is called
doubly-even if all of its words are annihilated by the first and the
second elementary symmetric polynomials in $n$ variables. For $q=2$ this
condition is actually equivalent to the usual one on weights modulo 4, but
for $q \geq 4$ it does not restrict the Hamming weight over $\F_q$. (For
odd $q$ the condition just means that the code is self-orthogonal and its
dual contains the all-ones word; however here we consider only
characteristic 2.) Extended Reed-Solomon codes of rate
$\frac{1}{2}$ are known to be examples of doubly-even self-dual codes.
For $q=4^e$ another interesting class of examples
is given by the extended quadratic-residue codes of lengths divisible
by 4.

We find (Theorems \ref{inv} and \ref{struct}) that the complete weight enumerators of
doubly-even self-dual codes over $\F_q$, $q=2^f$, are invariants for the
same type of Clifford-Weil group that for odd primes $q$ has been
discussed in \cite{RS98}, Section 7.9.
More precisely, the group has a normal subgroup of order
$4q^2$ or $8q^2$ (depending on whether $f$ is even or odd) such that
the quotient is $\SL_2(\F_q)$. Over $\F_4$ the invariant ring is still
simple enough to be described explicitly. Namely, the subring 
of Frobenius-invariant elements is
generated by the algebraically independent weight
enumerators of the four extended quadratic-residue codes of lengths 4, 8, 12
and 20, and the complete invariant ring is a free module of rank 2 over
this subring; the fifth (not Frobenius-invariant) basic generator has
degree 40. In the final section we use this result to find the maximal
Hamming distance of doubly-even self-dual quaternary codes up through length 24.
Over the field $\F _4$, doubly-even codes coincide with what are called
``Type II'' codes in \cite{Gab}.

The invariant ring considered here is always generated by weight
enumerators. This property holds even for
Clifford-Weil groups associated with multiple weight enumerators, for
which a direct proof in the binary case was given in \cite{cliff1}. The
general case can be found in \cite{cliff2}, where still more general
types of codes are also included.

\section{Doubly even codes}
In this section we generalize the notion of doubly-even
binary codes to arbitrary finite fields of characteristic 2
(see \cite{quebbemann}).

Let $\F := \F _{2^f}$ denote the field with $2^f$ elements.
A {\em code} $C \leq \F^{\,n}$ is an $\F$-linear subspace of $\F ^{\,n}$.
If $c \in \F ^{\,n}$ then the $i$-th coordinate of $c$ is denoted by $c_i$.
The {\em dual code}
to a code $C \leq \F ^{\,n}$
is defined to be
$$C^{\perp } := \{ v \in \F ^{\,n} \mid \sum _{i=1}^n c_iv_i =0 \fa
c \in C \} .$$
$C$ is called {\em self-orthogonal} if $C \subset C^{\perp}$, and {\em self-dual}
if $C=C^{\perp}$.

\begin{defn}
A code $C \leq \F ^{\,n} $ is {\em doubly-even}
if $$\sum _{i=1}^n c_i =
\sum _{i<j} c_ic_j = 0 $$
for all $c \in C $.
\end{defn}

\begin{rem}
An alternative definition can be obtained as follows.  There is a
unique unramified extension $\hat\F$ of the $2$-adic integers with the
property that $\hat\F/2\hat\F\cong\F$; moreover, the map $x\mapsto x^2$
induces a well-defined map $\hat\F/2\hat\F\to \hat\F/4\hat\F$, and thus
a map (also written as $x\mapsto x^2$) from $\F\to\hat\F/4\hat\F$.  The
above condition is then equivalent to requiring that $\sum_i v_i^2=0\in
\hat\F/4\hat\F$ for all $v\in C$.
\end{rem}

Doubly-even codes are self-orthogonal. This follows from the identity
$$\sum_{i<j} (c_i + c_i') (c_j + c_j') =
\sum_{i<j} c_i c_j + \sum_{i<j} c_i' c_j' + \sum_{i=1}^n c_i \sum_{i=1}^n c_i'
- \sum_{i=1}^n c_i c_i'.$$

Note that Hamming distances in a doubly-even code are not necessarily even:

\begin{exmp}
Let $\omega \in \F_4 $ be a primitive cube root of unity.
Then the code
$Q_4 \leq \F _4^{\,4}$ with generator matrix
$$ \left( \ba{cccc} 1  & 1 & 1 & 1 \\ 0 & 1 & \omega & \omega^2 \ea \right) $$
is a doubly-even self-dual code over $\F _4$.
\end{exmp}

Let $B=(b_1,\ldots , b_f )$ be an $\F _2$-basis of $\F$ such that
$\trace (b_ib_j) = \delta _{ij} \fa i,j=1,...,f$,
where $\trace $ denotes the trace of $\F$ over $\F_2$.
Then $B$ is called a
{\em self-complementary}
(or
{\em trace-orthogonal})
basis of $\F $
(cf. \cite{Pasq81},  \cite{quebbemann}, \cite{Wolf85}).
Using such a basis we identify
$\F  $ with $\F _2 ^{\,f}$ and define $$\varphi : \F \rightarrow \Z/4\Z, \
\varphi (\sum _{i=1}^f a_i b_i ) := \wt (a_1,\ldots , a_f ) + 4\Z $$
to be the weight modulo 4.
Since $\trace (b_i) = \trace (b_i^2) = 1$, we have
$$\varphi (a) + 2\Z = \trace(a)$$ and (considering $2\tau$ as a map onto
$2\Z/4\Z$)
$$ \varphi(a+a') = \varphi(a) + \varphi(a') + 2 \trace(aa')$$
for all $a,a' \in \F$. More generally,
$$
\varphi( \sum_ {i=1}^n c_i) = \sum _{i=1}^n \varphi(c_i) + 
2 \trace(\sum_{i<j} c_i c_j) .
$$
We extend $\varphi $ to a quadratic function
$$\phi : \F ^{\,n} \rightarrow \Z / 4\Z , \ \phi(c) := \sum _{i=1}^n \varphi (c_i). $$

\begin{prop}{\label{Phi}}
A code $C \leq \F^{\,n}$ is doubly-even if and only if $\phi (C) = \{0\} $.
\end{prop}

\proof
For $r\in \F, c \in \F^{\,n}$,
$$\phi(rc) = \varphi (\sum _{i=1}^n rc_i) - 2\trace (\sum _{i<j} r^2 c_i c_j ) .$$
This equation in particular shows that $\phi(C)=\{0\}$ if $C$ is doubly-even.
Conversely, if $\phi (C) = \{ 0 \} $ then the same equation shows that
$\trace(r \sum_{i=1}^n c_i) =\varphi (\sum _{i=1}^n rc_i) + 2\Z =0$
for all $r\in \F$, $c\in C$.
Since the trace bilinear form is non-degenerate, this implies that
$\sum _{i=1}^n c_i = 0$ for all $c\in C$.
The same equality then implies that
$\trace (r^2 \sum _{i<j} c_i c_j ) = 0 $ for all $r\in \F$ and $c\in C$.
The mapping $r \mapsto r^2$ is an automorphism of $\F$, so again
the non-degeneracy of the trace bilinear form yields $\sum _{i<j} c_i c_j = 0$
for all $c\in C$.
\eb

\begin{cor}
Let $\F ^{\,n}$ be identified with $\F _2^{\,nf}$ via a self-complementary
basis. Then a doubly-even code $C\leq \F^{\,n}$ becomes a doubly-even
binary code $C_{\F_2}\leq \F_2^{\,nf}$.
\end{cor}

\begin{rem}
Let $C\leq \F ^{\,n} $ be a doubly-even code.
Then $\1 := (1,\ldots ,1) \in C^{\perp }$.
Hence if $C$ is self-dual then $4$ divides $n$.
\end{rem}

In the following remark we use the fact that the length of a doubly-even self-dual
binary code is divisible by 8.

\begin{rem}
If $f \equiv 1 \pmod{2}$ then the length of a doubly-even self-dual code over $\F$
is divisible by $8$.
If $f \equiv 0 \pmod{2}$ then $\F \otimes _{\F _4 } Q_4$
is a doubly-even self-dual code over $\F$ of length $4$.
\end{rem}

More general examples of doubly-even self-dual codes are provided by 
extended quadratic-residue codes (see \cite{MS77}).
Let $p$ be an odd prime and let $\zeta $ be a primitive $p$-th
root of unity in an extension field $\tilde{\F} $ of $\F _2$.
Let
$$g:= \prod_{a\in (\F_p^*)^2} ( X - \zeta ^a ) \in \tilde{\F} [X ] $$
where $a$ runs through the non-zero squares in $\F_p $.
Then $g\in \F_4[X]$ divides $X^p-1$, and $g$ lies in $\F_2[X] $
if $g$ is fixed under the Frobenius automorphism $z\mapsto z^2$,
i.e. if $2$ is a square in $\F_p^*$, or
equivalently by quadratic reciprocity if $p\equiv \pm 1 \pmod{8} $.
Assuming $f$ to be even if $p\equiv \pm 3 \pmod{8}$,
 we define the {\em quadratic-residue code}
$\QR (\F,p) \leq \F ^{\,p}$ to be the cyclic code of length $p$
with generator polynomial $g$.
Then $\dim (\QR (\F , p )) = p - \deg (g) = \frac{p+1}{2} $,
which is also the dimension of
 the {\em extended code} $\widetilde{\QR} (\F , p) \leq \F ^{\,p+1}$.

   From \cite[pages 490, 508]{MS77}
together with Proposition \ref{Phi}
we obtain the following (the case $\F = \F _4$ was given in 
\cite[Proposition 4.1]{Gab}):

\begin{prop}{\label{QR}}
Let $p$ be a prime, $p\equiv 3 \pmod{4}$.
Then the extended quadratic-residue code
$\widetilde{\QR} (\F , p)$ is a doubly-even self-dual code.
\end{prop}

\section{Complete weight enumerators and invariant rings}

In this section we define the action of a group of $\C$-algebra
automorphisms on the polynomial ring $\C[x_a \mid a \in \F]$ such that
the complete weight enumerators of doubly-even self-dual codes are 
invariant under this group.

\begin{defn}
Let $C\leq \F ^{\,n}$ be a code. Then
$$\cwe (C) := \sum _{c \in C} \prod _{i=1}^n x_{c_i} \in \C [x_a \mid a \in \F ]$$
is the {\em complete weight enumerator} of $C$.
\end{defn}

For an element $r \in \F$ let $m_r$ and $d_r$ be the $\C$-algebra endomorphisms
of $ \C [x_a \mid a \in \F ]$ defined by
$$m_r(x_a) := x _{ar}, \ \
d_r(x_a) := i^{\varphi(ar)} x _{a} \ \
\fa a \in \F, $$
where $i= \sqrt{-1}$ and $\varphi : \F \rightarrow \Z / 4\Z $ is defined as above via
a fixed self-complementary basis.
We also have the MacWilliams transformation $h$ defined by
$$h(x_a) :=  2^{-f/2} \sum _{b\in \F } (-1) ^{\trace (ab) } x _b
\fa a \in \F .$$

\begin{defn}
The group $$G_f := \langle h, m_r , d_r \mid 0 \neq r\in \F \rangle $$
is called the associated {\em Clifford-Weil group}.
\end{defn}

Gleason (\cite{Gle70}) observed that the complete weight enumerator of a
self-dual code $C$ remains invariant under the transformations $h$ and $m_r$.
If $C$ is doubly-even, then $\cwe(C)$ is invariant also under each $d_r$
(Proposition \ref{Phi}). Therefore we have the following theorem.

\begin{thm}{\label{inv}}
The complete weight enumerator of a doubly-even self-dual code over $\F$
lies in the invariant ring
$$\Inv (G_f) := \{p \in
\C [x_a \mid a \in \F ] \mid p g = p \fa g \in G_f \}.$$
\end{thm}

By the general theory developed in \cite{cliff2} one finds that
a converse to Theorem \ref{inv} also holds:

\begin{thm}{\label{all}}
The invariant ring of $G_f$ is generated by complete weight
enumerators of doubly-even self-dual codes over $\F $.
\end{thm}

In the case $f=1$ Gleason obtained the more precise information
$$\Inv (G_1) = \C [ \cwe(\mathcal{H}_8), \cwe (\mathcal{G}_{24} ) ] $$
where $\mathcal{H}_8$ and $\mathcal{G}_{24}$
denote the extended Hamming code of length 8 and the
extended Golay code of length 24 over $\F_2$.

In general, the Galois group $$\Gamma_f := \Gal (\F / \F_2 )$$ acts on $\Inv(G_f)$
by $\gamma (x_a) := x_{a^\gamma} $ for all $a\in \F , \gamma \in \Gamma_f $.
Let $\Inv(G_f,\Gamma_f )$ denote the ring of $\Gamma_f $-invariant polynomials
in $\Inv(G_f)$.

\begin{thm}{\label{invf4}}
$$\Inv(G_2,\Gamma_2 ) = \C [ \cwe(Q_4) , \cwe (Q_8) ,
\cwe (Q_{12}), \cwe (Q_{20}) ] $$
where $Q_{p+1}$ denotes the extended quadratic-residue code of length
$p+1$ over $\F _4$ (see Proposition \ref{QR}).
The invariant ring of $G_2$ is a free module of rank $2$ over
$\Inv(G_2,\Gamma_2)$:
$$\Inv(G_2) =
\Inv(G_2,\Gamma_2 )  \oplus
\Inv(G_2,\Gamma_2 ) p_{40} $$
where $p_{40}$ is a homogeneous polynomial of degree $40$ which is not
invariant under $\Gamma _2$.
\end{thm}

\proof
Computation shows that $\langle G_2, \Gamma _2 \rangle$ is a complex reflection
group of order $2^9 3 \cdot 5$
(Number 29 in \cite{ShephardTodd}) and $G_2$ is a subgroup of index 2 with
Molien series
$$\frac{1+t^{40}}{(1-t^4)(1-t^8)(1-t^{12})(1-t^{20})} .$$
By Proposition \ref{QR} the codes $Q_i$ ($i=4,8,12,20$)  are 
 doubly-even self-dual codes over $\F_4$.
Their complete weight enumerators (which are $\Gamma _2$-invariant)
 are algebraically independent elements in the
invariant ring of $G_2$ as one shows by an explicit computation of
their Jacobi matrix.
Therefore these polynomials generate the algebra $\Inv (G_2, \Gamma_2)$.
\eb

By Theorem \ref{all} we have the following corollary.
\begin{cor}
There is a doubly-even self-dual code $C$ over $\F_4$ of length $40$
such that $\cwe(C)$ is not Galois invariant.
\end{cor}

A code with this property was recently constructed in \cite{BeCh02}.

For $f>2$ the following example shows that we cannot
hope to find an explicit description of the invariant
rings of the above type.

\begin{exmp}
The Molien series of $G_3$ is
$N/D$, where
$$ D = (1-t^8)^2(1-t^{16})^2(1-t^{24})^2(1-t^{56})(1-t^{72})$$
and $N(t) = M(t) + M(t^{-1})t^{216}$ with
$$M =
1+
 5t^{16} +
 77t^{24} +
 300t^{32} +
 908t^{40} +
 2139t^{48} +
 3808t^{56}+
 5864t^{64}$$ $$+
 8257t^{72} +
     10456t^{80} +
 12504t^{88} +
 14294t^{96} +
 15115t^{104}.$$
The Molien series of $\langle G_3, \Gamma_3 \rangle$ is $(L(t)+L(t^{-1})t^{216})/D$,
where $D$ is as above and
$$L=
 1 + 3t^{16 } + 29t^{24 } + 100t^{32 } + 298t^{40 } + 707t^{48 } +
 1268t^{56 }+ 1958t^{64 }$$ $$
 + 2753t^{72 }
   + 3482t^{80 }+ 4166t^{88 } + 4766t^{96 } + 5045t^{104 }.$$
\end{exmp}

\section{The structure of the Clifford-Weil groups $G_f$}

In this section we establish the following theorem.

\begin{thm}{\label{struct}}
The structure of the Clifford-Weil groups $G_f$ is given by
$$G_f \cong Z.(\F \oplus \F ).\SL_2(\F ) $$
where $Z \cong \Z/4\Z $ if $f$ is even, and
$Z \cong \Z/8\Z $ if $f$ is odd.
\end{thm}

To prove this theorem, we first construct a normal subgroup
$N_f \unlhd G_f$ with $N_f\cong \Z/4\Z \zentr 2_{+}^{1+2f} $, the central
product of an extraspecial group of order $2^{1+2f}$ with the cyclic group of 
order 4.
The image of the homomorphism $G_f/N_f \rightarrow \Out (N_f)$  is isomorphic
to $\SL _2 (\F )$ and the kernel consists of scalar matrices only.

Let $q_r := (d_r^2)^h = hd_r^2h$ and
$$N_f:= \langle d_r^2 , q_r , i \id \mid r\in \F \rangle .$$
Using the fact that
$(-1)^{\varphi (b)} = (-1)^{\trace(b)} $ for all $b\in \F$, we
find that
$$d_r^2(x_a) = (-1)^{\tau(ar)}x_a, \ \ q_r(x_a) = x_{a+r} .$$
For the chosen self-complementary basis
 $(b_1,\ldots , b_f)$,
 $q_{b_j}$ commutes with $d_{b_k}^2$ if $j\neq k$ and the commutator of
$q_{b_j}$ and $d_{b_j}^2$ is $-\id $.  From this we have:

\begin{rem}{\label{irred}}
The group $N_f$ is isomorphic to a central product of an extraspecial group
$\langle q_{b_j} , d_{b_j}^2 \mid j=1,\ldots , f \rangle \cong 2^{1+2f}_+$
with the center $Z(N_f) \cong \Z/4\Z$. The representation of $N_f$
on the vector space $\oplus _{a\in \F } \C x_a $ of dimension $2^f$ is
the unique irreducible
representation of $N_f$ such that $t\in \Z/4\Z$ acts as multiplication
by $i^t$.
\end{rem}

Concerning the action of $G_f$ on $N_f$ we have
$$ m_a  d_r^2  m_a^{-1} = d_{a^{-1} r}^2 , \ \
 m_a  q_r  m_a^{-1} = q_{a r} , \ \fa a,r \in \F^{\,*} .$$
Since $m_a$ conjugates $d_r$ to $d_{a^{-1}r}$ it suffices to calculate
the action of $d_1$:
$$ d_1  d_r^2  d_1^{-1} = d_{ r}^2 , \ \
 d_1  q_r  d_1^{-1} = i^{\varphi(r)} q_r d_r^2  , \ \fa r \in \F .$$
This proves

\begin{lem}
The image of the homomorphism $G_f \rightarrow \Aut (N_f/Z(N_f)) $
is isomorphic to $\SL_2(\F )$ via
$$h\mapsto \left(\ba{cc} 0 & 1 \\ 1 & 0 \ea \right),\ \
m_a \mapsto \left(\ba{cc} a & 0 \\ 0 &  a ^{-1} \ea \right),\ \
d_1 \mapsto \left( \ba{cc} 1 & 0 \\ 1 & 1 \ea \right) .$$
\end{lem}

Elementary calculations or explicit knowledge of the automorphism group of $N_f$
(see \cite{Win72}) show that the kernel
of the above homomorphism is $N_f C_{G_f} (N_f) = N_f (G_f \cap \C^* \id ) $.
It remains to find the center of $G_f$, which
by the calculations above contains $i \id $.
If $f$ is even, then  $\cwe(Q_4 \otimes _{\F_4 } \F )$ is an invariant of
degree 4 of $G_f$, so the center is isomorphic to $\Z/4\Z$ in this case.
To prove the theorem,
it remains to construct an element $\zeta _8 \id \in G_f $ if $f$ is odd, where
$\zeta _8 \in \C^*$ is a primitive $8$-th root of unity.

\begin{lem}
If $f$ is odd, then
$\langle (hd_1 )^3 \rangle = \langle \zeta _8 \id \rangle $.
\end{lem}

\proof
$(hd_1)^3$ acts trivially on $N_f/Z(N_f) $.
Explicit calculation shows that $(hd_1)^3$ commutes with each generator
of $N_f$, hence acts as a scalar.
We find that
$$(hd_1)^3 (x_0) = \frac{1}{\sqrt{|\F |}} \frac{1}{|\F| }\sum _{b,c\in \F} i^{\varphi (c+b)} (-1) ^{\trace (c) } x_0 .$$
The right hand side is an $8$-th root of unity times $x_0$.
If $f$ is odd, then $\sqrt{2} $ is mentioned, which implies that this
is a primitive $8$-th root of unity.
\eb

\section{Extremal codes}

Let $C \leq \F ^{\,n}$ be a code. The complete weight enumerator
$\cwe(C) \in \C [ x_a \mid a \in \F ]$
may be used to obtain information about the Hamming weight enumerator,
which is the polynomial in a single variable $x$ obtained from $\cwe(C)$ by
substituting
$x_0 \mapsto 1$ and $x_a \mapsto x \fa a \neq 0.$

\begin{rem}
(a)
If $\F ' \leq \F $ is a subfield of $\F$ and $e = [\F  :\F']$, then
$C$ becomes a code $C_{\F '}$ of length $en$ over $\F '$ by
identifying $\F $ with $\F'^{\,e}$ with respect to a self-complementary basis
$(b_1,\ldots , b_e)$.
If $a = \sum _{i=1}^e a_i b_i$ with $a_i\in \F'$, then the complete
weight enumerator of $C_{\F'}$ is obtained from $\cwe(C)$ by
replacing $x_a$ by $\prod_{i=1}^e x_{a_i} $.
\\
(b)
We may also construct a code $C'$ of length $n$ over $\F '$ from
$C$ by taking the $\F' $-rational points:
$$C':= \{ c\in C \mid c_i \in \F' \fa i=1,\ldots , n \}.$$
The dimension of $C'$ is at most the dimension of $C$, and
the complete weight enumerator of $C'$ is found by the substitution
$x_a \mapsto 0$ if $a\not\in \F '$.
$C'$ is called the {\em $\F'$-rational subcode} of $C$.
\end{rem}

As an application of Theorem \ref{invf4} we have the following result.
Note that the results for lengths $n \le 20$ also follow from the classification
of doubly-even self-dual codes in \cite{Gab}, \cite{BetGuHM01}
and \cite{Bets02}, and the bound 
for length 20 can be deduced from \cite[Cor. 3.4]{Gab}.

\begin{thm}
Let $\F := \F_4$.
The maximal Hamming distance $d=d(C)$ of a doubly-even self-dual code
$C\leq \F^{\,n}$ is as given in the following table:
$$\begin{array}{|c|c|c|c|c|c|c|}
\hline
n & 4 & 8 & 12 & 16 & 20 & 24  \\
\hline
d & 3 & 4 & 6 & 6 & 8  & 8 \\ \hline \end{array}
$$
For $n=4$ and $8$, the quadratic-residue codes $Q_4$ resp. $Q_8$ are the
unique codes $C$ of length $n$ with $d(C) = 3$ resp. $d(C) = 4$.
\end{thm}

\proof
Let $p \in \C [x_0,x_1,x_{\omega } , x_{\omega ^2}]_{n}^{G_2} $,
a homogeneous polynomial of degree $n$.
If $p$ is the complete weight enumerator of a code $C$ with $d(C) \geq d$, then
the following conditions must be satisfied.
\begin{itemize}
\item[a)] All coefficients in $p$ are nonnegative integers.
\item[b)] The coefficients of $x_0^ax_1^bx_{\omega }^b x_{\omega  ^2}^b$ with
$b>0$ are divisible by 3.
\item[c)] $p(1,1,1,1) = 2^n $.
\item[d)] $p(1,1,0,0) = 2^m$ for some $m\leq \frac{n}{2} $.
\item[e)] $p(1,x,x,x) - 1 $ is divisible by $x^d$.
\end{itemize}
One easily sees that $Q_4$ is the unique doubly-even self-dual code over
$\F$ of length 4.
If $C$ is such a code of length 8 with $d(C) \geq 4$, then
$\cwe(C) $ is uniquely determined by condition e).
In particular the $\F_2$-rational subcode of $C$ has dimension 4 and is
a doubly-even self-dual binary code of length 8. Hence
$C= \mathcal{H}_8 \otimes \F = Q_8$.
If $C \leq \F^{\,12} $ is a doubly-even self-dual code with $d(C) \geq 6$,
then again $\cwe(C) = \cwe (Q_{12})$ is uniquely determined by condition
e), moreover $Q_{12} $ has minimal distance $6$.

For $n=16$, there is a unique polynomial
$p(x_0,x_1,x_{\omega } , x_{\omega^2} ) \in \C [x_0,x_1,x_{\omega } , x_{\omega ^2}]_{16}^{G_2} $
such that $p(1,x,x,x) \equiv 1+ax^7 \pmod{x^8}$.
This polynomial $p$ has negative coefficients.
Therefore the doubly-even self-dual codes $C\leq \F^{\,16}$ satisfy
$d(C) \leq 6$. There are two candidates for polynomials
$p$ satisfying the five conditions above with $d=6$.
The rational subcode has either dimension 2 or 4 and all words
$\neq 0, {\bf 1}$ are of weight 8.
One easily constructs such a code $C$ from the code $Q_{20}$,
by taking those elements of $Q_{20}$ that have 0 in four fixed coordinates,
omitting these 4 coordinates to get a code of length 16, adjoining the all-ones
vector and then a vector of the form $(1^8,0^8)$ from the dual code.
$C_{\F_2} \leq \F_2^{\,32}$ is isomorphic to the extended binary 
quadratic-residue code and the rational subcode of $C$ is 2-dimensional.

For $n=20$ we similarly find four candidates for complete weight enumerators
satisfying a) - e) above with $d=8$ (where the
dimension of the rational subcode is $1,3,5$ or $7$).
None of these satisfies e) with $d>8$.
 The code $Q_{20}$ has minimal weight 8 and
its rational subcode is $\{ {\bf 0} , {\bf 1} \} $.
For $n=24$, the code $Q_{24}
= \F_4 \otimes \mathcal{G}_{24}$ has $d(C) = 8$.
To see that this is best possible let
$p\in \C[x_0,x_1,x_{\omega }, x_{\omega ^2}]^{G_2}_{24}$
satisfy (b) and (e) above
with $d=9$.
Then $p=p_0 + ah_1+ b h_2$, for suitable  $p_0,h_1,h_2$
with $h_i(1,x,x,x) \equiv 0 \pmod{x^9}$,
$p_0(1,x,x,x) \equiv 1 \pmod{x^9}$ and $a,b \in \Z $.
Explicit calculations then show that
 $p_0(1,1,0,0)$, $h_1(1,1,0,0)$ and $h_2(1,1,0,0)$ are all divisible by 3.
Therefore $p(1,1,0,0)$ is not a power of 2, hence $p$ does not satisfy condition d).
\eb

\vspace{1.0\baselineskip}
{\bf Acknowledgment.} We thank O. Jahn for computations in
connection with Theorem \ref{invf4}  at an early stage of this work.
We also thank the referees for their comments.


\begin{thebibliography}{1234}

\bibitem{Bets02}
K. Betsumiya,
On the classification of Type II codes over $\F _{2^r}$
with binary length 32, preprint, 2002.

\bibitem{BeCh02}
K. Betsumiya and Y. J. Choie,
Codes over $\F_4$, Jacobi forms and Hilbert-Siegel
modular forms over $\Q(\sqrt{5})$, preprint, 2002.

\bibitem{BetGuHM01}
K. Betsumiya, T. A. Gulliver, M. Harada and A. Munemasa,
On type II codes over $\F_4$, {\it IEEE Trans. Inform. Theory} {\bf 47},
(2001), 2242--2248.

\bibitem{Gab} 
P. Gaborit, V. S. Pless, P. Sol\'e and A. O. L. Atkin, 
Type II codes over $\F _4$,
{\it Finite Fields Applic.} {\bf 8} (2002), 171-183.

\bibitem{Gle70}
A. M. Gleason,
Weight polynomials of self-dual codes and the MacWilliams identities,
in {\it Actes, Congr\'{e}s International de Math\'{e}matiques (Nice, 1970)},
Gauthiers-Villars,
Paris, 1971, Vol. 3, pp. 211--215.

\bibitem{klemm}
M. Klemm, Eine Invarianzgruppe f\"{u}r die vollst\"{a}ndige Gewichtsfunktion
selbstdualer Codes, {\it Archiv Math.} {\bf 53} (1989), 332-336.

\bibitem{MS77}
F. J. MacWilliams and N. J. A. Sloane,
{\it The Theory of Error-Correcting Codes},
North-Holland, Amsterdam, 1977.

\bibitem{cliff1}
G. Nebe, E. M. Rains and N. J. A. Sloane,
The invariants of the Clifford groups,
{\it Designs, Codes Cryptography} {\bf 24} (2001), 99--121.

\bibitem{cliff2}
G. Nebe, E. M. Rains and N. J. A. Sloane,
{\it Self-Dual Codes and Invariant Theory},
book in preparation.

\bibitem{Pasq81}
G. Pasquier,
Binary self-dual codes construction from
self-dual codes over a Galois field $\F _{2^m}$,
in {\it Combinatorial Mathematics (Luminy, 1981)},
ed. C. Berge et al., Annals Discrete Math. {\bf 17} (1983), 519--526.

\bibitem{quebbemann}
H.-G. Quebbemann, On even codes,
{\it Discrete Math.}  {\bf 98}  (1991), 29--34.


\bibitem{RS98}
E. M. Rains and N. J. A. Sloane,
{\it Self-dual codes}
in {\it Handbook of Coding Theory},
ed. V. Pless and W. C. Huffman,
Elsevier, Amsterdam, 1998, pp. 177--294.


\bibitem{ShephardTodd} G. C. Shephard and J. A. Todd,
 Finite unitary reflection groups,
{\it Canad. J. Math.} {\bf 6} (1954), 274-304.

\bibitem{Win72}
D. L. Winter,
The automorphism group of an extraspecial $p$-group,
{\it Rocky Mtn. J. Math.} {\bf 2} (1972), 159--168.

\bibitem{Wolf85}
J. Wolfmann,
A class of doubly even self-dual binary codes,
{\it Discrete Math.} {\bf 56} (1985), 299--303.

\end{thebibliography}
\end{document}